\documentclass{amsart}
\usepackage{amssymb,amsmath,amsfonts,amscd,euscript}

\newcommand{\C}{\mathbb C}

\newcommand{\N}{\mathbb N}
\newcommand{\Z}{\mathbb Z}
\newcommand{\Pro}{\mathbb P}

\newcommand{\T}{\otimes}
\newcommand{\ld}{\ldots}

\newcommand{\hk}{\hookrightarrow}
\newcommand{\hki}{\hookleftarrow}

\newcommand{\nc}{\newcommand}

\nc{\slt}{\mathfrak{sl}_2}
\nc{\g}{\mathfrak{g}}
\nc{\slth}{\widehat{\slt}}
\nc{\Slth}{\widehat{\Slt}}

\nc{\al}{\alpha}
\nc{\be}{\beta}
\nc{\veps}{\varepsilon}

\nc{\ov}{\overline}

\nc{\codim}{{\mathop{\rm codim}}}
\nc{\Id}{{\mathop{\rm Id}}}
\nc{\im}{{\mathop{\rm im}}}
\nc{\const}{{\mathop{\rm const}}}
\nc{\sltc}{\slt\otimes\C[t]}
\nc{\bra}{\langle}
\nc{\ket}{\rangle}
\nc{\Slt}{{\mathop{\rm SL}}_2}
\nc{\G}{{\mathop{\rm G}}}

\nc{\Sltc}{\Slt (\C[t]/t^n)}
\nc{\sh}{{\mathop{\rm sh}}}
\nc{\Pic}{{\mathop{\rm Pic}}}
\nc{\U}{{\mathop{\rm U}}}
\nc{\Fl}{{\mathop{\rm Fl}}}
\nc{\is}{\{i_1,\ld,i_s\}}
\nc{\an}{a_1,\ld,a_n}
\nc{\Sum}{\sum\nolimits}

\nc{\shn}{\sh^{(n)}}
\nc{\shin}{\sh_i^{(n-1)}}

\nc{\Ob}{\EuScript {O}}
\nc{\E}{\EuScript {E}}

\newtheorem{rem}{Remark}[section]
\newtheorem*{rem*}{Remark}
\newtheorem{opr}{Definition}[section]
\newtheorem{theorem}{Theorem}[section]
\newtheorem*{theorem*}{Theorem}
\newtheorem{utv}{Statement}[section]
\newtheorem{lem}{Lemma}[section]
\newtheorem{cor}{Corollary}[section]
\newtheorem{prop}{Proposition}[section]

\begin{document}
\pagestyle{plain}
\markboth{E.Feigin}
{Schubert varieties and the fusion products.\\ The general case.}
\title[Schubert varieties and the fusion products]
{Schubert varieties and the fusion products.\\ The general case.}
\author{E.Feigin}
\address{Independent University of Moscow, Russia, Moscow,
Bol'shoi Vlas'evskii per.,7. and \\ 
Moscow State University, Russia, Moscow, Leninskie gory, 1,Faculty of 
Mathematics and Mechanics, Department of the Higher Algebra}
\email {evgfeig@mccme.ru }
\date{}
\begin{abstract}
This paper generalizes the results of the paper \cite{mi3} to the case of
the general $\slt$ Schubert varieties. We study the homomorphisms between
different
Schubert varieties, describe their geometry and the group of the line
bundles. We also derive some consequences concerning the infinite-dimensional
generalized affine grassmanians.
\end{abstract}

\maketitle
\section*{Introduction}
Let $A=(a_1,\ld,a_n)\in\N^n, 1<a_1\le\ld\le a_n$. Let $M^A$ be the fusion
product of
$\slt$ modules $\C^{a_1},\ld,\C^{a_n}$ (see \cite{fusion, mi1}) (recall that
$M^A$ is $\slt\T (\C[t]/t^n)$ module). Let $v_A$ be the lowest weight
vector in $M^A$ with respect to the $h_0$-grading (for $x\in\slt$ denote
$x_i=x\T t^i$). Denote by $[v_A]$ the line
$\C\cdot [v_A]\in\Pro(M^A)$.
In \cite{mi3} the closure $\sh_A=\ov{\Slt(\C[t]/t^n)\cdot [v_A]}\hk
\Pro(M^A)$
was studied in the case when $a_i\ne a_j$ for $i\ne j$.
It was proved that in this case the Schubert variety $\sh_A$ is smooth
$n$-dimensional algebraic variety,
independent on the choice of $A$ (the only demand is $a_i\ne a_j$).
This variety was denoted by $\shn$.
In the present paper we study the case of the general Schubert variety.

We start with the generalization of the independence of
$\sh_A$ ($a_i\ne a_j$)
on $A$. We say that $A$ is of the type $\is$ if
$i_1+\ld +i_s=n$ and
$$a_1=\ld =a_{i_1}\ne a_{i_1+1}=\ld =a_{i_1+i_2}\ne\ld\ne
a_{i_1+\ld+i_{s-1}+1}=\ld =a_n.$$
Then $\sh_A\simeq \sh_B$ if and only if $A$ and $B$ are of the same
type. We denote the corresponding variety by $\sh_{\is}$ and the point
$[v_A]$ by $[v_{\is}]$. For example,
$\shn=\sh_{\{1,\ld,1\}}$.

In \cite{mi3} for any $n>k$ the bundle $\shn\to \sh^{(k)}$ with a fiber
$\sh^{(n-k)}$ was constructed. The generalization of this construction is the
following fact: let $1\le t<s$. Then there exists
$\Slt(\C[t]/t^n)$-equivariant bundle $\sh_{\is}\to\sh_{\{i_{t+1},\ld,i_s\}}$
with a fiber $\sh_{\{i_1,\ld,i_t\}}$, sending $[v_{\is}]$ to
$[v_{\{i_{t+1},\ld,i_s\}}]$. For the proof we study some special
subvarieties of $\sh_{\is}$. Namely, we prove that
$$\sh_{\is}=\Sltc\cdot [v_{\is}]\cup\bigcup_{j=1}^{s-1}
N_{i_1+\ld +i_j}(\is),$$
where $N_\al(\is)$ are some subvarieties of $\sh_{\is}$. The latter are
closely connected
with $\slt\T (\C[t]/t^n)$ submodules $S_{i,i+1}(A)\hk M^A$,
studied in \cite{mi2, mi3}.
Recall that these submodules were defined via the following exact sequence of
$\slt\T (\C[t]/t^n)$ modules (see \cite{SW1, SW2} for the similar relations on the
$q$-characters):
$$0\to S_{i,i+1}(A)\to M^A\to  M^{(a_1,\ld,a_i-1, a_{i+1}+1,\ld,a_n)}\to 0.$$
Now let
$$A_{\is}=(2^{i_1}\ld (s+1)^{i_s})=(\underbrace{2,\ld,2}_{i_1},\ld,
\underbrace{(s+1),\ld,(s+1)}_{i_s}).$$
Fix an isomorphism $\sh_{A_{\is}}\simeq \sh_{\is}$.
Then $$N_\al(\is)\hk \Pro (S_{\al,\al+1}(A_{\is})).$$

An important property is the existence of the surjective
homomorphisms between the different Schubert varieties. Namely, we say that
$\is\ge \{j_1,\ld,j_{s_1}\}$ if there exists such numbers
$1\le k_1<\ld < k_{s_1-1}< s$ that
$$i_1+\ld +i_{k_1}=j_1,\ld, i_{k_{s_1-1}+1}+\ld +i_s=j_{s_1}.$$
For example, $\{1,\ld,1\}$ is the largest type, while $\{n\}$ is the
smallest one. We will prove that if
$\is\ge \{j_1,\ld,j_{s_1}\}$ then there exists a surjective
$\Slt(\C[t]/t^n)$-equivariant
homomorphism $\sh_{\is}\to \sh_{\{j_1,\ld,j_{s_1}\}}$.

Schubert varieties have a description in terms of some partial flags in
$W_0=\C^2\T\C[t]$ (for the analogous construction in the case of $\shn$
see \cite{mi3}; see also \cite{segal}).
$W_0$ is naturally $\sltc$ module and we also have
an action of the operator $t$ by multiplication. Consider the
variety $\Fl_{\is}$ of the sequences of the subspaces of $W_0$:
\begin{multline*}
\Fl_{\is}=\{W_0\hki W_1\hki\ld\hki W_s:\\ 1).\ tW_\al\hk W_\al;\quad
2).\ \dim W_\al/W_{\al+1}=i_{s-\al};\quad
3).\ W_{\al+1}\hki t^{i_{s-\al}}W_\al\}.
\end{multline*}
We prove that $\Fl_{\is}\simeq\sh_{\is}$.

We study the Picard group of the Schubert varieties. Let $\E$ be the line
bundle on $\sh_{\is}$. Following \cite{mi3} denote by $C_i$ the projective
lines:
$$C_i=\ov{\left\{\exp (ze_i)\cdot [v_{\is}], z\in\C\right\}}\hk
\sh_{\is}, i=0,\ld,n-1.$$
One can show that $\E$ is
completely determined by its restriction to these lines.
Let $B=(b_1,\ld,b_n)\in\Z^n$.
Introduce a notation
$\E=\Ob(B)$ if $\E|_{C_i}\simeq \Ob(b_1+\ld+b_{n-i})$.
We prove that the bundle $\Ob(B)$ really exists on $\sh_{\is}$
if and only if the type of $B$ is less or equal to $\is$.
It is possible to describe the space of
sections of $\Ob(B)$ in terms of the fusion products. Namely
if $b_i\ge 0$ then
$$H^0(\Ob(B),\sh_{\is})\simeq
\left(M^{(b_1+1,\ld,b_n+1)}\right)^*.$$
As a consequence
we obtain the proof of the theorem about the sections of the line bundles on
the generalized
affine grassmanians (see \cite{mi3, sto, Gordon}).
Recall the main definitions. An embedding
$\sh_{\is}\hk \sh_{\{i_1,\ld,i_s+2\}}$ allows us to organize an inductive
limit
$$Gr_{\is}=\lim_{k\to\infty}
\sh_{(i_1,\ld,i_s+2k)}.$$
There exists a bundle $\Ob(B^{(\infty)})$
on $Gr_{\is}$ such that
$$\Ob(B^{(\infty)})|_{\sh_{\{i_1,\ld,i_s+2k\}}}=
\Ob(b_1,\ld,b_n,\underbrace{b_n,\ld,b_n}_{2k}).$$
We decompose
the space of sections of the above bundle as $\slth$ module
in the case $0\le b_1\le\ld\le b_n$:
$$H^0(\Ob(B^{(\infty)}), Gr_{\is})\simeq
\bigoplus_{j=0}^{b_n} c_{j;b_1,\ld,b_n} L^*_{j,b_n}.$$
Here $L_{j,b_n}$ are irreducible $\slth$ modules and $c_{j;b_1,\ld,b_n}$
-- the structure constants of the level $b_n+1$ Verlinde algebra,
associated with the Lie algebra $\slt$.

We finish the paper with the discussion of the singularities of the Schubert
varieties. Namely we prove that the only smooth Schubert variety is $\shn$
and study the singularities of the "smallest" variety $\sh_{\{n\}}$.

The paper is organized in the following way:

Section $1$ contains the preliminary statements from the papers
\cite{mi1, mi2, mi3} and some generalizations (lemma
$(\ref{S_i}))$.

In the section $2$ we identify the isomorphic Schubert varieties
(theorem $(\ref{isom})$).

Section $3$ is devoted to the proof of the existence of the bundle
$\sh_{\is}\to \sh_{\{i_{t+1},\ld,i_s\}}$ with a fiber $\sh_{\{i_1,\ld,i_t\}}$
(theorem $(\ref{mainth})$).

In the section $4$ we
give the description of $\sh_{\is}$ in terms of the generalized partial
flag manifolds (proposition $(\ref{flag})$).

Section  $5$ is devoted to the study of the line bundles on $\sh_{\is}$
(proposition $(\ref{linebun})$) and
the spaces of their sections (corollary $(\ref{fus})$).

In the section $6$ we decompose the spaces of sections of some line bundles
on the
infinite-dimensional generalized affine grassmanians into the sum of the
dual irreducible $\slth$ modules (proposition $(\ref{inf})$).

Section $7$ contains the discussion of the singularities of the Schubert
varieties.

{\bf Acknowledgment.} I want to thank B.Feigin for useful discussions.
This work was partially supported by the grants RFBR  03-01-00167 and
SS 1910.2003.1.

\section {Preliminaries}
Here we briefly recall the main notions and statements about the Schubert
varieties $\sh_A$ from \cite{mi3}.

Let $A=(\an)\in\N^n, 1<a_1\le\ld\le a_n$, $M^A$ the corresponding fusion
product
(see \cite{fusion, mi1}), $v_A$ and $u_A$ its lowest and highest weight
vectors with respect to the
$h_0$-grading (for $x\in\slt$ $x_i=x\T t^i$). $M^A$  is cyclic
$\slt\T(\C[t]/t^n)$
module,
$M^A= \C[e_0,\ld,e_{n-1}]\cdot v_A.$ The group $\G=\Sltc$ acts on $M^A$ and
thus
on its projectivization $\Pro(M^A)$. Schubert variety $\sh_A\hk\Pro(M^A)$
is the closure of the orbit of the point $[v_A]$ (for $v\in M^A$
$[v]=\C\cdot v\in \Pro(M^A))$:
$$\sh_A=\ov{\G\cdot [v_A]}.$$
It was proved in \cite{mi3} that $\sh_A$ is projective complex algebraic
variety. Its coordinate ring can be described in the following way. Recall
that for $A, B, C\in\N^n$ with $c_i=a_i+b_i-1$ the multiplication
$(M^A)^*\T (M^B)^*\to (M^C)^*$ was constructed. Thus, for any $A$ we have
an algebra $$CR_A=\bigoplus_{i=0}^{\infty} (M^{A_i})^*,$$ where
$A_i=(ia_1-i+1,\ld, ia_n-i+1),\ i\ge 0.$
It was proved in \cite{mi3} that $CR_A$ is a coordinate ring of
$\sh_A$.

We will need the realization of the fusion products in the tensor powers
of the space of the semi-infinite forms (fermionic realization),
constructed in \cite{mi2}. Recall
that the space of the semi-infinite forms $F$ (the fermionic space)
is the level $1$ $\slth$ module.
$F$ contains the set of vectors $v(i), i\in\Z$ with the following properties:
\begin{gather*}
1.\ \ e_{i-1}v(i)=v(i-2),\qquad 2.\ \ e_kv(i)=0 \text{ if } k\ge i,\\
3.\ \ \U(\slth)\cdot v(0)\simeq L_{0,1},\qquad
4.\ \   \U(\slth)\cdot v(1)\simeq L_{1,1}.
\end{gather*}
It was proved in \cite{mi2} that $M^A$ can be embedded into
$F^{\T (a_n-1)}$. Namely,
\begin{equation}
\label{emb}
M^A\simeq \U(\sltc)\cdot \left(v(n)\T v(n-d_2)\T\ld\T
v(n-d_2-\ld -d_{a_n-1})\right),
\end{equation}
where $d_i=\#\{\alpha:\ a_{\alpha}=i\}$.

Recall the special submodules of $M^A$, studied in \cite{mi3}.
For any $i, 1\le i< n$ where exists an $\sltc$ submodule $S_{i,i+1}(A)\hk
M^A$ with a property
$$M^A/S_{i,i+1}(A)\simeq
M^{(a_1,\ld,a_{i-1},a_i-1,a_{i+1}+1, a_{i+2},\ld,a_n)}.$$
There are three cases when the modules $S_{i,i+1}(A)$ can be easily described:
\begin{gather*}
1.\ \ S_{1,2}(A)\simeq M^{(a_2-a_1+1,a_3,\ld,a_n)};\\
2.\ \  \text{ if } a_i=a_{i+1}, \text { then }
S_{i,i+1}(A)\simeq M^{(a_1,\ld,a_{i-1},a_{i+2},\ld,a_n)};\\
3.\ \ S_{n-1,n}(A)\simeq M^{(a_1,\ld,a_{n-2})}\otimes \C^{a_n-a_{n-1}+1}.
\end{gather*}
In the general case we have the following lemma, which follows
from the fermionic realization of the fusion product (for the analogous
proofs see \cite{mi3}).
\begin{lem}
\label{S_i}
Let $a_i<a_{i+1}$.
Denote
$$A'=(a_1,\ld,a_{i-1},a_{i+2},\ld,a_n),\qquad
A''=(a_{i+1}-a_i+1,\ld,a_n-a_i+1).$$
Then
there is an embedding of $\sltc$ modules
$$S_{i,i+1}(A)\hk M^{A'}\T M^{A''}.$$
The image of this embedding coincides with
\begin{equation}
\label{submod}
\C[e_0,\ld,e_{n-1},e^{(2)}_{n-i-1}]\cdot \left( v_{A'}\T v_{A''}\right),
\end{equation}
where $e_j^{(i)}$ stands for the operator $e_j$ acting on the $i$-th
factor of the tensor product. From now on we identify $S_{i,i+1}(A)$ with
its image $(\ref{submod})$. The vector $v_{A'}\T v_{A''}$ (as well as
its preimage in $S_{i,i+1}(A)$) is denoted by $v_{i,i+1}(A)$.
\end{lem}

\section{The main definition}
Let $A=(a_1\le\ld\le a_n)\in\Z^n$. We say that $A$ is of the type
$\{i_1,\ld,i_s\}$ with $i_1+\ld +i_s=n, i_\alpha\ge 1$ if
\begin{equation*}
a_1=\ld =a_{i_1}\ne a_{i_1+1}=a_{i_1+2=}\ld =a_{i_1+i_2}\ne\ld\ne
a_{n-i_s+1}=\ld =a_n.
\end{equation*}
Now let $A,B\in (\N\setminus 1)^n$.
We will prove that if $A$ and $B$ are of the same type, then
$\sh_A\simeq\sh_B$.
Introduce a notation: if $A$ is of the type $\{i_1,\ld,i_s\}$, then
we denote
$$A=(a_1^{i_1}a_{i_1+1}^{i_2}\ld a_{i_1+\ld +i_{s-1}+1}^{i_s}).$$

\begin{lem}
\label{lemisom}
$\sh_{(2^n)}\simeq\sh_{(k^n)}$ for any $k\ge 2$.
\end{lem}
\begin{proof}
Recall that the coordinate ring of $\sh_{(2^n)}$ is
$\bigoplus_{i=0}^{\infty} (M^{(i^n)})^*$. The coordinate ring of
$\sh_{(k^n)}$ is $\bigoplus_{j=0}^{\infty} (M^{((jk-j+1)^n)})^*$. Thus
$\sh_{(2^n)}\simeq\sh_{(k^n)}$ (see \cite{alggeom}).
\end{proof}

\begin{theorem}
\label{isom}
Let $A,B\in (\N\setminus 1)^n$ be of the same type. Then $\sh_A\simeq \sh_B$.
\end{theorem}
\begin{proof}
It is enough to show that if $A$ is of the type $\{i_1,\ld,i_s\}$, then
$$\sh_A\simeq \sh_{(2^{i_1}3^{i_2}\ld (s+1)^{i_s})}.$$
Note that there is a natural surjective $\G$-equivariant homomorphism
$$\sh_A\to \sh_{(2^{i_1}3^{i_2}\ld (s+1)^{i_s})},\quad
[v_A]\mapsto [v_{(2^{i_1}3^{i_2}\ld (s+1)^{i_s})}].$$
For the proof we use the
fermionic realization of the fusion product. Let $\phi_A$ be the embedding
$M^A\hk F^{\T (a_n-1)}$ (sometimes we write simply $\phi$ instead of
$\phi_A$).
Then, because of $(\ref{emb})$
\begin{multline}
\label{emb1}
\phi_A(v_A)=
v(n)^{\T (a_1-1)}\T\\ \T v(n-i_1)^{\T (a_{i_1+1}-a_{i_1})}\T\ld\T
v(n-i_1-\ld-i_{s-1})^{\T(a_{n-i_s+1}-a_{i_s})}.
\end{multline}
Picking one factor from each tensor power
$$v(n-i_1-\ld -i_\alpha)^{\T
(a_{i_1+\ld +i_\alpha+1}-a_{i_1+\ld +i_\alpha})}$$
we obtain a vector
$$\phi_{(2^{i_1}3^{i_2}\ld (s+1)^{i_s})}(v_{(2^{i_1}3^{i_2}\ld (s+1)^{i_s})})
\in F^{\T s}.$$
Thus for some permutation $\sigma\in S_{a_n-1}$ of the factors of
$F^{\T (a_n-1)}$ and a
vector $w\in F^{\T (a_n-s-1)}$ we have
\begin{equation*}
\sigma\phi_A(v_A)=
\phi_{(2^{i_1}3^{i_2}\ld (s+1)^{i_s})}(v_{(2^{i_1}3^{i_2}\ld (s+1)^{i_s})})
\T w.
\end{equation*}
Hence because of $\sh_A=\ov{\G\cdot [v_A]}$ we obtain a surjective map
$$\sh_A\to \sh_{(2^{i_1}3^{i_2}\ld (s+1)^{i_s})},\quad
[v_A]\mapsto [v_{(2^{i_1}3^{i_2}\ld (s+1)^{i_s})}].$$
We need to prove that this map is an isomorphism.

 Note that from $(\ref{emb1})$ we obtain an embedding
$$\sh_A\hk \sh_{((a_1)^n)}\times \sh_{((a_{i_1+1}-a_{i_1}+1)^{n-i_1})}
\times\ld
\times \sh_{((a_n-a_{n-i_s}+1)^{n-i_1-\ld -i_{s-1}})}.$$
Also we have an embedding
$$\sh_{(2^{i_1}3^{i_2}\ld (s+1)^{i_s})}\hk
\sh_{(2^n)}\times \sh_{(2^{n-i_1})}\times\ld\times
\sh_{(2^{n-i_1-\ld -i_{s-1}})}$$
and the following diagram is commutative:
$$
\begin{CD}
\sh_A  @>>>  \sh_{((a_1)^n)}\times\ld
\times \sh_{((a_n-a_{n-i_s}+1)^{n-i_1-\ld -i_{s-1}})} \\
@VVV  @V{\psi}VV  \\
\sh_{(2^{i_1}3^{i_2}\ld (s+1)^{i_s})} @>>>
\sh_{(2^n)}\times \sh_{(2^{n-i_1})}\times\ld\times
\sh_{(2^{n-i_1-\ld -i_{s-1}})}.
\end{CD}
$$
Note that because of the lemma $(\ref{lemisom})$, the map $\psi$ is an
isomorphism as the product of the isomorphisms
$$\sh_{((a_{i_1+\ld+i_\alpha+1}-a_{i_1+\ld+i_\alpha}+1)^{n-i_1-\ld-i_\alpha})}
\simeq \sh_{(2^{n-i_1-\ld-i_\alpha})}.$$
Thus the left vertical map from the diagram is an isomorphism.
\end{proof}

\begin{opr}  Let $A\in (\N\setminus 1)^n$ be of the type $\{i_1,\ld,i_s\}$.
Define  $\sh_{\{i_1,\ld,i_s\}}=\sh_A$. Denote
$[v_A]=[v_{\{i_1,\ld,i_s\}}],\ \ [u_A]=[u_{\{i_1,\ld,i_s\}}]$.
\end{opr}

\noindent {\bf Example.}  Recall that in \cite{mi3} the case of the
$\sh_A$ with pairwise distinct $a_i$ was considered. It was proved that
all the Schubert varieties of this type are isomorphic. This variety was
denoted as $\shn$. In our notations
$\shn=\sh_{\{\underbrace{\scriptstyle{1,\ld,1}}_n\}}$.
\newcounter{a}

\section{The existence of the bundle
$\sh_{\{i_1,\ld,i_s\}}\to \sh_{\{i_{t+1},\ld,i_s\}}$}
Recall that in \cite{mi3} for any
$n>k$ a bundle $\pi_{n,k}:\shn\to \sh^{(k)}$ with a
fiber
$\sh^{(n-k)}$ was constructed. In this section
we prove that for any $t<s$ there exists $\G$-equivariant bundle
$$\sh_{\{i_1,\ld,i_s\}}\to \sh_{\{i_{t+1},\ld,i_s\}}\quad \text
{ with a fiber } \quad
\sh_{\{i_1,\ld,i_t\}}.$$ In order to do that, we study the structure of
some special subvarieties of $\sh_{\{i_1,\ld,i_s\}}$.

\begin{lem}
\label{glue}
Let $1\le \alpha <s$. Then there is a $\G$-equivariant surjective
homomorphism
$$\sh_{\is}\to \sh_{\{i_1,\ld, i_{\alpha-1}, i_\alpha+i_{\alpha+1},
i_{\alpha+2},\ld,i_s\}}, \quad
[v_{\is}]\mapsto [v_{\{i_1,\ld, i_\alpha+i_{\alpha+1}, \ld,i_s\}}],
$$
and the restriction
$G\cdot [v_{\is}]\to \G\cdot
[v_{\{i_1,\ld, i_\alpha+i_{\alpha+1}, \ld,i_s\}}]$ is an isomorphism.
\end{lem}

\begin{proof}
Fix the realizations
\begin{equation*}
\sh_{\is}=\sh_{(2^{i_1}\ld (s+1)^{i_s})},\quad
\sh_{\{i_1,\ld,i_\al+i_{\al+1},\ld,i_s\}}=
\sh_{(2^{i_1}\ld (\al+1)^{i_\al+i_{\al+1}}\ld s^{i_s})}.
\end{equation*}
Note that
\begin{multline}
\label{vec1}
\phi(v_{(2^{i_1}\ld (s+1)^{i_s})})=
v(n)\T v(n-i_1)\T\ld\T v(n-i_1-\ld-i_{s-1});\\
\phi (v_{(2^{i_1}\ld (\alpha+1)^{i_\alpha+i_{\alpha+1}}\ld s^{i_s})})
=v(n)\T \ld \T v(n-i_1-\ld -i_{\alpha-1})\T \\ \T
v(n-i_1-\ld-i_\alpha-i_{\alpha+1})\T\ld\T  v(n-i_1-\ld-i_{s-1}).
\end{multline}
For the sequence $0< n_1<\ld < n_j\le n$ define a homomorphism of
$\slth$ modules
$P(n_1,\ld,n_j): F^{\T n}\to F^{\T j}$:
$$
P(n_1,\ld,n_j) (w_1\T\ld\T w_n) =w_{n_1}\T\ld \T w_{n_j}.
$$
From the formula $(\ref{vec1})$ we obtain
\begin{equation}
\label{P}
P(1,2,\ld,\alpha, \alpha+2,\ld,s) \phi(v_{(2^{i_1}\ld (s+1)^{i_s})})=
\phi (v_{(2^{i_1}\ld (\alpha+1)^{i_\alpha+i_{\alpha+1}}\ld s^{i_s})}).
\end{equation}
The formula $(\ref{P})$ gives us the needed homomorphism.
\end{proof}

Introduce a notation: let $i_1+\ld +i_s=j_1+\ld+j_{s_1}=n$. We write
$\is\ge \{j_1,\ld,j_{s_1}\}$ if for any $A, B$ of the types $\is$,
$\{j_1,\ld,j_{s_1}\}$ correspondingly the following condition holds:
$(a_i=a_{i+1}\Rightarrow b_i=b_{i+1})$.
For example, $\{1,\ld,1\}$ is a maximal type, while $\{n\}$ is a minimal one.

\begin{cor}
\label{his}
Let $\is\ge \{j_1,\ld,j_{s_1}\}$. Then there exists a $\G$-equivariant
surjective homomorphism
$\sh_{\is}\to \sh_{\{j_1,\ld,j_{s_1}\}}$.
In particular, for any type $\is$ there exists a homomorphism
$h_{\is}:\shn\to \sh_{\is}$.
\end{cor}

For the proof of the existence of the bundle $\sh_{\is}\to
\sh_{\{i_{t+1},\ld,i_s\}}$ we need to study the structure of the complement
$\sh_{\is}\setminus \G\cdot [v_{\is}]$.
First recall some results from \cite{mi3}.
Denote $B_n=(2,\ld,n+1)$. Thus
$B_n$ is of the type $\{1,\ld,1\}$ and $\sh_{B_n}= \shn$.
\begin{utv}
There exist $(n-1)$-dimensional subvarieties $N_j(\{1,\ld,1\})\hk\shn$,
$j=1,\ld,n-1$ with the following properties:\\
$1.\
\shn\setminus \G\cdot [v_{B_n}]=  \bigcup_{j=1}^{n-1}
N_j (\{1,\ld,1\});$\\
$2.\
N_j(\{1,\ld,1\})\simeq  \{(x,y)\in\sh^{(n-2)}\times\sh^{(n-j)}:\
\pi_{n-2,n-j-1}(x)=\pi_{n-j,n-j-1}(y)\};$\\
$3.\
N_j(\{1,\ld,1\})\hk \Pro(S_{j,j+1}(B_n));$\\
$4.\
N_j(\{1,\ld,1\})=
\ov{\{\exp \left(\Sum_{i=0}^{n-1} e_iz_i +e^{(2)}_{n-j-1}z\right)\cdot
[v_{j,j+1}(B_n)], z_i,z\in\C\}}$ (the closure is taken in the
projective space $\Pro(S_{j,j+1}(B_n))$; about the notation
$e^{(2)}_{n-j-1}$ see lemma $(\ref{S_i})$).
\end{utv}

We will need the following lemma.
\begin{lem}
\label{homlem}
Denote $A_{\is}=(2^{i_1}\ld (s+1)^{i_s})$. Then
\begin{equation*}
h_{\is} [v_{j,j+1}(B_n)]=
[v_{i_1+\ld +i_\alpha, i_1+\ld +i_\alpha+1}(A_{\is})]
\end{equation*}
for all $j$ with $i_1+\ld +i_{\alpha-1}< j \le i_1+\ld +i_\alpha$
(we put $i_0=0$).
\end{lem}
\begin{proof}
Recall (see \cite{mi3}) that
\begin{multline*}
\phi_{B_n} (v_{j,j+1}(B_n))=\\ =v(n-2)\T v(n-3)\T \ld \T v(n-j-1)\T v(n-j)\T
v(n-j-1) \T\ld \T v(0).
\end{multline*}
At the same time
\begin{multline*}
\phi_{A_{\is}} (v_{i_1+\ld +i_\alpha, i_1+\ld +i_\alpha+1}(A_{\is}))=
v(n-2)\T v(n-i_1-2)\T\ld \\ \ld\T v(n-i_1-\ld -i_{\alpha-1}-2)\T
v(n-i_1-\ld -i_\alpha)\T\ld \T v(n-i_1-\ld-i_{s-1}).
\end{multline*}
Thus for $i_1+\ld +i_{\alpha-1}< j \le i_1+\ld +i_\alpha$
\begin{multline*}
P(1,i_1+1,i_1+i_2+1,\ld,i_1+\ld+ i_{s-1}+1) (\phi_{B_n} (v_{j,j+1}(B_n)))= \\
=\phi_{A_{\is}} (v_{i_1+\ld +i_\alpha, i_1+\ld +i_\alpha+1}(A_{\is})).
\end{multline*}
This gives us
$h_{\is} [v_{j,j+1}(B_n)]=
[v_{i_1+\ld +i_\alpha, i_1+\ld +i_\alpha+1}(A_{\is})]$.
\end{proof}

Now let $\is$ be some type. Fix a realization
$\sh_{\is}=\sh_{A_{\is}}$. Define subvarieties of $\sh_{\is}$:
\begin{gather*}
N_j(\is)=
\ov{\left\{\exp\left(\sum\nolimits_{l=0}^{n-1} z_le_l +
e^{(2)}_{n-j-1}z\right)\cdot
[v_{j,j+1}(A_{\is})], z_l,z\in\C\right\}},\\
j=i_1,i_1+i_2,\ld,i_1+\ld +i_{s-1}.
\end{gather*}
From the lemma $(\ref{homlem})$ we obtain the following corollary:
\begin{cor}
Let $s>1$. Then
$$\sh_{\is}\setminus \G\cdot [v_{\is}]= \bigcup_{\alpha=1}^{s-1}
N_{i_1+\ld +i_\alpha} (\is).$$
\end{cor}
\begin{proof}
Recall the surjective homomorphism $h_{\is}: \shn\to \sh_{\is}$. Note that
the restriction
$$h_{\is}: \G\cdot [v_{B_n}] \to \G\cdot [v_{\is}]$$
is an isomorphism. In addition, lemma $(\ref{homlem})$ gives us that
\begin{gather*}
h_{\is} N_j(\{1,\ld,1\})=N_{i_1+\ld +i_\alpha}(\is),\\
i_1+\ld +i_{\alpha-1}< j \le i_1+\ld +i_\alpha.
\end{gather*}
Corollary is proved.
\end{proof}

For the description of the varieties $N_j(\is)$  we need the following notation:
let $\al_1,\ld,\al_{s_1},\be_1,\ld,\be_{s_2}$ be the natural numbers
with
$$\al_1+\ld+\al_{s_1}=n-2,\quad \be_1+\ld +\be_{s_2}=n-j.$$
Consider the map
\begin{equation*}
h_{\{\al_1,\ld,\al_{s_1}\}}\times h_{\{\be_1,\ld,\be_{s_2}\}}:
\sh^{(n-2)}\times\sh^{(n-j)}\to
\sh_{\{\al_1,\ld,\al_{s_1}\}}\times \sh_{\{\be_1,\ld,\be_{s_2}\}}.
\end{equation*}
Recall that $N_j(\{1,\ld,1\})\hk \sh^{(n-2)}\times\sh^{(n-j)}$.
Define the variety
\begin{equation*}
\sh_{\{\al_1,\ld,\al_{s_1}\}}\ \widetilde{\times}\
\sh_{\{\be_1,\ld,\be_{s_2}\}}=
h_{\{\al_1,\ld,\al_{s_1}\}}\times h_{\{\be_1,\ld,\be_{s_2}\}}
(N_j(\{1,\ld,1\})).
\end{equation*}
The following proposition gives the description of the varieties
$N_j(\is)$.

\begin{prop}
\label{mainprop}
Let $A=A_{\is}$.
\begin{enumerate}
\item  $s=1$. Then $\sh_{\{n\}}=\G\cdot [v_A]\sqcup \sh_{\{n-2\}}$.
\item  $s>1$. Then
$\sh_{\is}= \G\cdot [v_A]\sqcup \bigcup_{\alpha=1}^{s-1}
N_{i_1+\ld +i_\alpha} (\is)$ and
\begin{enumerate}
\item\label{>2} $i_\al \ge 2$. Then
\begin{equation*}
N_{i_1+\ld +i_\alpha} (\is)\simeq
\sh_{\{i_1,\ld,i_\al-2,\ld,i_s\}}\ \widetilde{\times}\
\sh_{\{i_{\al+1},\ld,i_s\}}
\end{equation*}
(if $i_\al=2$ then $\{i_1,\ld,i_\al-2,\ld,i_s\}=\{i_1,\ld,i_{\al-1},
i_{\al+1},\ld,i_s\}$).
\item\label{=1} $i_\al=1$. Then
\begin{equation*}
N_{i_1+\ld +i_\alpha} (\is)\simeq
\sh_{\{i_1,\ld,i_{\al-1},i_{\al+1}-1,\ld,i_s\}}\ \widetilde{\times}\
\sh_{\{i_{\al+1},\ld,i_s\}}.
\end{equation*}
\end{enumerate}
\end{enumerate}
\end{prop}
\begin{proof}
Let $s=1$. Then $A=(2^n)$. Consider the map $h_{\{n\}}:\shn\to \sh_{\{n\}}$
and the vector $v_A'=e_{n-1} v_A\in M^A$. Then for any $j=1,\ld,n-1$
$h_{\{n\}} (v_{j,j+1}(B_n))={v'_A}$. The proof is similar to the one from
the lemma $(\ref{homlem})$. Thus
$$h_{\{n\}} N_j(\{1,\ld,1\})=\ov{G\cdot [{v'}_A]}.$$
But
$\phi_A ({v'}_A)=v(n-2)$. Hence $\ov{G\cdot [{v'}_A]}\simeq \sh_{\{n-2\}}$
and
$$\sh_{\{n\}}=\G\cdot [v_A]\sqcup \ov{G\cdot [{v'}_A]}\simeq
\G\cdot [v_A]\sqcup \sh_{\{n-2\}}.$$

Now let $s>1$ and $i_\al\ge 2$.
Because of the lemma $(\ref{S_i})$ and the definition of the varieties
$N_j(A)$ we obtain the embedding
\begin{equation*}
N_{i_1+\ld +i_\alpha} (\is)\hk
\sh_{\{i_1,\ld,i_\al-2,\ld,i_s\}} \times
\sh_{\{i_{\al+1},\ld,i_s\}}.
\end{equation*}
At the same time the following diagram is commutative
$$
\begin{CD}
N_{i_1+\ld+i_\al}(\{1,\ld,1\}) @>>>
\sh^{(n-2)}\times \sh^{(n-i_1-\ld-i_\al)}\\
@V{h_{\is}}VV @VV{h_{\{i_1,\ld,i_\al-2,\ld,i_s\}}\times h_{\{i_{\al+1},\ld,i_s\}}}V\\
N_{i_1+\ld+i_\al}(A) @>>>
\sh_{\{i_1,\ld,i_\al-2,\ld,i_s\}} \times
\sh_{\{i_{\al+1},\ld,i_s\}}.
\end{CD}
$$
This finishes the proof of the case $(\ref{>2})$. The proof of the case
$(\ref{=1})$ is quite similar.
\end{proof}

As a corollary, we prove the theorem about the existence of the bundles:
\begin{theorem}
\label{mainth}
Fix a type $\is$, $i_1+\ld +i_s=n$. For any \  $t=1,\ld,s-1$ there is a
$G$-equivariant bundle
$\sh_{\is}\to \sh_{\{i_{t+1},\ld,i_s\}}$ with a fiber
$\sh_{\{i_1,\ld,i_t\}}$, sending $[v_{\is}]$ to $[v_{\{i_{t+1},\ld,i_s\}}]$.
\end{theorem}
\begin{proof}
Note that there is a natural $\G$-equi\-va\-ri\-ant surjective homomorphism
$$\Psi_t: \sh_{\is}\to \sh_{\{i_{t+1},\ld,i_s\}},\quad
\Psi_t([v_{\is}])=[v_{\{i_{t+1},\ld,i_s\}}],
$$
because
\begin{multline*}
\phi_{A_{\is}} (v_{A_\is})=\\ =
v(n)\T v(n-i_1)\T\ld \T v(n-i_1-\ld-i_t)\T
\phi_{A_{\{i_{t+1},\ld,i_s\}}} (v_{A_{\{i_{t+1},\ld,i_s\}}}).
\end{multline*}
We want to prove that for any $x\in \sh_{\{i_{t+1},\ld,i_s\}}$ the preimage
$\Psi_t^{-1}(x)$ is isomorphic to $\sh_{\{i_1,\ld,i_t\}}$.

First prove that
$\Psi_t^{-1}([v_{\{i_{t+1}\ld i_s\}}])\simeq \sh_{\{i_1,\ld,i_t\}}$.
Recall (see \cite{mi1}) that
$$\C[e_{i_{t+1}+\ld+i_s},\ld,e_{n-1}]\cdot v_{\is}\simeq
M^{(2^{i_1}\ld (t+1)^{i_t})}.$$
Thus we obtain
\begin{multline*}
\Psi_t^{-1}([v_{\{i_{t+1}\ld i_s\}}])=
\ov{\left\{\exp\left(\sum\nolimits_{j=i_{t+1}+\ld +i_s}^{n-1} e_jz_j\right)
\cdot [v_{\is}], z_j\in\C\right\}} \simeq \\ \simeq
\ov{\left\{\exp\left(\sum\nolimits_{j=0}^{i_1+\ld+i_t-1} e_jz_j\right)\cdot
[v_{\{i_1,\ld,i_t\}}], z_j\in\C\right\}}\simeq \sh_{\{i_1,\ld,i_t\}}.
\end{multline*}
In the same way one can prove that the preimage $\Psi_t^{-1}(x)$ is
isomorphic to $\sh_{\{i_1,\ld,i_t\}}$ for any $x$ from the orbit
$\G\cdot [v_{\{i_{t+1},\ld,i_s\}}]$.

Recall that
$$\sh_{\{i_{t+1},\ld,i_s\}}= \G\cdot [v_{\{i_{t+1},\ld,i_s\}}]\bigsqcup
\bigcup_{\al=1}^{s-t-1} N_{i_{t+1}+\ld +i_{t+\al}}(\{i_{t+1},\ld,i_s\}).$$
Pick $\al$ with $1\le \al \le s-t-1$. Let $i_{t+\al}\ge 2$.
Note that
$$\Psi_t^{-1}:\sh_{\{i_{t+1},\ld,i_\al-2,\ld,i_s\}}\ \widetilde{\times} \
\sh_{\{i_{\al+1},\ld,i_s\}}=
\sh_{\{i_1,\ld,i_{\al}-2,\ld,i_s\}}\ \widetilde{\times} \
\sh_{\{i_{\al+1},\ld,i_s\}}
$$
and for any
$$x\in \sh_{\{i_1,\ld,i_{\al}-2,\ld,i_s\}},\quad
y\in \sh_{\{i_{\al+1},\ld,i_s\}}$$ we have
$\Psi_t (x\times y)=\Psi_t(x) \times y$.
The assumption that
$$\Psi_t: \sh_{\{i_1,\ld,i_{\al}-2,\ld,i_s\}}\to
\sh_{\{i_{t+1},\ld,i_\al-2,\ld,i_s\}}$$
is a bundle with a fiber $\sh_{\{i_1,\ld,i_t\}}$ gives us
that for any
$$x\in N_{i_{t+1}+\ld +i_{t+\al}}(\{i_{t+1},\ld,i_s\})
\text{ with } i_{t+\al}\ge 2$$
the preimage $\Psi_t^{-1} (x)$ is isomorphic to
$\sh_{\{i_1,\ld,i_t\}}$. The case of $i_{t+\al}=1$ can be considered
in the same way. Theorem is proved.
\end{proof}

\begin{cor}
Let $P_{\is}(q)=\sum_{j=0}^{2n} q^j\dim H_j(\sh_{\is},\Z)$. Then
$$P_{\is}(q)=\prod_{\al=1}^s \frac{(1-q^{2(i_\al+1)})}{(1-q^2)}.$$
\end{cor}
\begin{proof}
Our corollary follows from the previous theorem and the statement that
\begin{equation}
\label{hom}
P_{\{n\}}(q)=\frac{1-q^{2(n+1)}}{1-q^2},
\end{equation}
i.e.
the odd homologies vanish and $\dim H_{2j} (\sh_{\{n\}},\Z)=1,
\ j=0,\ld,n$. Because of the proposition $(\ref{mainprop})$
\begin{equation}
\label{dec}
\sh_{\{n\}}=\G\cdot [v_{\{n\}}]\sqcup\sh_{\{n-2\}}.
\end{equation}
But in
\cite{mi3} was proved that for any $A\in (\N\setminus)^n$ the orbit
$\G\cdot [v_A]\hk \sh_A$ is fibered over $\Pro^1$ with a fiber $\C^{n-1}$.
Thus $\G\cdot [v_{\{n\}}]$ is the union of the two cells: $\C^n$ and
$\C^{n-1}$. Using $(\ref{dec})$ we obtain $(\ref{hom})$.
\end{proof}

\begin{rem}
Consider the bundle $\G\cdot [v_A]\to \Pro^1$. It was proved in \cite{mi3}
that the transition functions of this bundle can be written in the
following form. Let $p\in\Pro^1\setminus \{0,\infty\}$, $x_0$ the coordinate
of $p$ in $\Pro^1\setminus \{\infty\}$, $y_0=x_0^{-1}$ the coordinate
of $p$ in $\Pro^1\setminus \{0\}$. Denote by  $\{x_i\}, \{y_j\},$
$i,j=1,\ld,n-1$ the
coordinates in the fiber $\C^{n-1}$ over the point $p$ with respect to
the trivializations on $\Pro^1\setminus \{\infty\}$ and $\Pro^1\setminus
\{0\}$. Then
we have an equality in the ring $(\C[t]/t^n)$:
$$(x_0+x_1t+\ld +x_{n-1}t^{n-1})(y_0+y_1t+\ld +y_{n-1}t^{n-1})=1.$$
The Schubert varieties
$\sh_{\is}$ can be considered as the $2^{n-1}$ ways of the com\-pactification
of the bundle over $\Pro^1$ with the fiber $\C^{n-1}$ and above transition
functions.
\end{rem}

\section {Schubert varieties as a generalized partial flag manifolds}
Here we give a description of the varieties $\sh_{\is}$ in terms of the
special sequences of the subspaces of $\C[t]\oplus \C[t]$.
We start with the case of $\sh_{\{n\}}$.
\begin{lem}
\label{s=1}
Let $W_0=\C^2\T\C[t]$ with a natural action of the operator $t$ by
multiplication. Define $\Fl_{\{n\}}$ as a variety of the
subspaces $W_1\hk W_0$ with the following properties:
\begin{equation}
\label{cond}
1).\ \dim W_0/W_1=n\qquad 2).\  tW_1\hk W_1\qquad 3).\ W_1\hki t^n W_0.
\end{equation}
Then $\Fl_{\{n\}}\simeq \sh_{\{n\}}$.
\end{lem}
\begin{proof}
First note that because of the conditions $1)$ and $3)$ we can consider
$W_1$ as a
subspace of $\C^2\T (\C[t]/t^n)$ of codimension $n$. Thus the group
$\G=\Sltc$ naturally acts on $\Fl_{\{n\}}$.
Let $v, u$ be the standard basis of the $2$-dimensional $\slt$ module,
$hv=-v,\ hu=u,\ ev=u.$ Denote $v_i=v\T t^i, u_i=u\T t^i, u_i,v_i\in W_0$.
Let
$V_{\{n\}}\in \Fl_{\{n\}}$ be the subspace with a basis $v_i, i=0,\ld,n-1$.
One can show that
\begin{equation}
\label{clos}
\Fl_{\{n\}}=\ov{\G\cdot V_{\{n\}}}.
\end{equation}
Consider a map
\begin{equation}
\label{wedge}
\Fl_{\{n\}}\to \bigwedge\nolimits^n \left(\C^2\T (\C[t]/t^n)\right),\quad
W_1\mapsto r_1\wedge\ld\wedge r_n,
\end{equation}
where $r_i$ is a basis of $W_1$. Surely $(\ref{wedge})$ is a
$\G$-equivariant homomorphism. Thus, because of the formula $(\ref{clos})$
for the proof of the lemma it is enough to show that
we have an isomorphism of $\slt\T (\C[t]/t^n)$ modules
\begin{equation}
\label{real}
\slt\T (\C[t]/t^n)\cdot (v_0\wedge\ld\wedge v_{n-1})\simeq
M^{(2^n)}.
\end{equation}
First we check that the defining relations from the right hand side of
$(\ref{real})$ are true in the left hand side. Recall (see \cite{mi1}) that
$$M^{(2^n)}\simeq \C[e_0,\ld,e_{n-1}]/I,$$
where $I$ is the ideal, generated by the following conditions
\begin{equation}
\label{eq}
e^{(n)}(z)^i=(e_{n-1}+ze_{n-2}+\ld +z^{n-1}e_0)^i\div z^{n(i-1)}.
\end{equation}
(The latter means that the first $n(i-1)-1$ coefficients of the series
$e^{(n)}(z)^i$ vanish.)
We will prove that
$e^{(n)}(z)^i\cdot (v_0\wedge\ld\wedge v_{n-1}) \div z^{n(i-1)}.$  (One can
easily check that the left hand side of $(\ref{real})$ will not change
after the replacement of $\slt\T (\C[t]/t^n)$ by $\C[e_0,\ld,e_{n-1}]$.)

By definition
\begin{multline}
\label{sum}
e^{(n)}(z)^i (v_0\wedge\ld\wedge v_{n-1})=\\
=n!\sum_{0\le \al_1<\ld <\al_i\le n-1} v_0\wedge\ld\wedge e^{(n)}(z)v_{\al_1}
\wedge\ld\wedge e^{(n)}(z) v_{\al_i}\wedge\ld\wedge v_{n-1}.
\end{multline}
Let us show that
$e^{(n)}(z)v_{\al_1}\wedge\ld\wedge e^{(n)}(z) v_{\al_i}\div z^{n(i-1)}$.
In fact
\begin{multline}
\label{form}
e^{(n)}(z)v_{\al_1}\wedge\ld\wedge e^{(n)}(z) v_{\al_i}=\\
=\sum_{j_1=\al_1}^{n-1} u_{j_1}z^{n+\al_1-j_1-1}\wedge
 \sum_{j_2=\al_2}^{n-1} u_{j_2}z^{n+\al_2-j_2-1}\wedge\ld
\wedge\sum_{j_i=\al_i}^{n-1} u_{j_i}z^{n+\al_i-j_i-1}.
\end{multline}
Note that
\begin{multline*}
\sum_{j_k=\al_k}^{n-1} u_{j_k}z^{n+\al_k-j_k-1}\wedge
\sum_{j_{k+1}=\al_{k+1}}^{n-1} u_{j_{k+1}}z^{n+\al_{k+1}-j_{k+1}-1}=\\
=\sum_{j_k=\al_k}^{\al_{k+1}-1} u_{j_k}z^{n+\al_k-j_k-1}\wedge
\sum_{j_{k+1}=\al_{k+1}}^{n-1} u_{j_{k+1}}z^{n+\al_{k+1}-j_{k+1}-1}.
\end{multline*}
Thus we can rewrite the formula $(\ref{form})$ in the following way:
\begin{multline*}
e^{(n)}(z)v_{\al_1}\wedge\ld\wedge e^{(n)}(z) v_{\al_i}=\\
=\sum_{j_1=\al_1}^{\al_2-1} u_{j_1}z^{n+\al_1-j_1-1}\wedge
 \sum_{j_2=\al_2}^{\al_3-1} u_{j_2}z^{n+\al_2-j_2-1}\wedge\ld
\wedge\sum_{j_i=\al_i}^{n-1} u_{j_i}z^{n+\al_i-j_i-1}\div\\
\div z^{n+\al_1-\al_2}z^{n+\al_2-\al_3}\ld z^{\al_i}=z^{n(i-1)+\al_1}\div
z^{n(i-1)}.
\end{multline*}

To finish the proof of the formula $(\ref{real})$ we show that the following
vectors are linearly independent ($k=0,\ld,n$):
$$e_{i_1}\ld e_{i_k} (v_0\wedge \ld \wedge v_{n-1}),\  i_\al\le n-k\  \
(1\le \al\le k).$$
(That will be enough, because $\dim M^{(2^n)}=2^n$.)
Let
\begin{equation}
\label{linind}
\sum_{0\le i_1\le\ld\le i_k\le n-k}\be_{i_1,\ld,i_k}
e_{i_1}\ld e_{i_k} (v_0\wedge\ld\wedge v_{n-1})=0.
\end{equation}
Pick such $n$-tuple $(i_1^0,\ld,i_k^0)$ that $\be_{i_1^0,\ld,i_k^0}\ne 0$ and
for any $(i_1,\ld,i_k)$ with a non-vanishing
$\be_{i_1,\ld,i_k}$ the following is true: there exists such $l\le k$
that $i_l>i^0_l$ and for any $m<l$ $i_m=i^0_m$. Then the monomial
$$e_{i_1^0} v_0\wedge\ld\wedge e_{i_k^0} v_{k-1}\wedge v_k\wedge\ld\wedge
v_{n-1}=
u_{i_1^0}\wedge u_{i_2^0+1}\wedge\ld\wedge u_{i_k^0+k-1}\wedge v_k\wedge
\ld\wedge v_{n-1}$$
comes from $e_{i_1^0}\ld e_{i_k^0} (v_0\wedge\ld\wedge v_{n-1})$ but
not from any other monomial from the linear combination $(\ref{linind})$.
Thus $(\ref{real})$ is proved.
\end{proof}

We will need the following lemma:
\begin{lem}
\label{embed}
Let $i_1+\ld +i_s=n$. Then there is $\G$-equivariant embedding
\begin{gather}
\sh_{\is}\hk \sh_{\{n\}}\times \sh_{\{n-i_1\}}\times\ld
\times \sh_{\{n-i_1-\ld-i_{s-1}\}},\\
\notag [v_{\is}]\mapsto [v_{\{n\}}]\times [v_{\{n-i_1\}}]\times\ld
\times [v_{\{n-i_1-\ld-i_{s-1}\}}].
\end{gather}
\end{lem}
\begin{proof}
Recall the embedding $M^{(2^{i_1}\ld (s+1)^{i_s})}\hk F^{\T s}$
$$v_{(2^{i_1}\ld (s+1)^{i_s})}\mapsto v(n)\T v(n-i_1)\T\ld \T
v(n-i_1-\ld-i_{s-1}).$$
Thus we obtain an embeddings
\begin{gather*}
M^{(2^{i_1}\ld (s+1)^{i_s})}\hk M^{(2^n)}\T M^{(2^{n-i_1})}\T\ld\T
M^{(2^{n-i_1-\ld-i_{s-1}})},\\
\sh_{\is}\hk \sh_{\{n\}}\times \sh_{\{n-i_1\}} \times\ld\times
\sh_{\{n-i_1-\ld-i_{s-1}\}}.
\end{gather*}
(Surely, the first embedding is a homomorphism of $\slt\T
(\C[t]/t^n)$ modules and the second one is a $\G$-equivariant homomorphism.)
\end{proof}

\begin{prop}
\label{flag}
Let $W_0=\C^2\T \C[t]$.
Define the generalized partial flag manifold $\Fl_{\is}$ as the variety of
the special sequences of the subspaces. Namely
\begin{multline*}
\Fl_{\is}=\{W_0\hki W_1\hki\ld\hki W_s:\\ 1).\ tW_\al\hk W_\al;\quad
2).\ \dim W_\al/W_{\al+1}=i_{s-\al};\quad
3).\ W_{\al+1}\hki t^{i_{s-\al}}W_\al\}.
\end{multline*}
Then $\Fl_{\is}\simeq\sh_{\is}$.
\end{prop}
\begin{proof}
Because of the condition $2)$ we know that $W_\al\hki t^n W_0$. Thus we can
consider $W_\al$ as a subspaces of $\C^2\T (\C[t]/t^n)$. Define
$V_{\is}\in \Fl_{\is}$ as follows:
\begin{gather*}
V_{\is}=\{W_0\hki W_1\hki\ld\hki W_s\},\\
W_\al=\bra v_i, i=0,\ld,n-1;\ u_j, j=i_s+\ld +i_{s-\al+1},\ld,n-1\ket,
1\le \al\le s-1,\\
W_s=\bra v_i, i=0,\ld,n-1\ket.
\end{gather*}
Note that
$\Fl_{\is}=\ov{\G\cdot V_{\is}}$. Define a map
\begin{multline}
\label{flmap}
\Fl_{\is}\to \bigwedge\nolimits^{n+i_1+\ld+i_{s-1}}
\left(\C^2\T (\C[t]/t^n)\right)\oplus\\ \oplus
\bigwedge\nolimits^{n+i_1+\ld+i_{s-2}}\left(\C^2\T (\C[t]/t^n)\right)\oplus
\ld\oplus \bigwedge\nolimits^{n}\left(\C^2\T (\C[t]/t^n)\right),\\
\{W_0\hki W_1\hki\ld\hki W_s\}\mapsto \bigoplus_{\al=1}^s
r^\al_1\wedge \ld \wedge r^\al_{n+i_1+\ld+i_{s-\al}},
\end{multline}
where $r^\al_j$ form a basis of $W_\al$.
Consider an embedding $\chi_\al, \al=1,\ld,s$:
\begin{gather*}
\chi_\al:\bigwedge\nolimits^{i_s+\ld +i_{s-\al+1}}
(\C^2\T (\C[t]/t^{i_s+\ld +i_{s-\al+1}}))\hk
\bigwedge\nolimits^{n+i_1+\ld + i_{s-\al}} (\C[t]/t^n),\\
\chi_\al(w)= w\wedge v_{i_s+\ld +i_{s-\al+1}}\wedge\ld\wedge v_{n-1}\wedge
u_{i_s+\ld +i_{s-\al+1}}\wedge\ld\wedge u_{n-1}.
\end{gather*}
One can show that the image of the map $(\ref{flmap})$ belongs to the direct
sum of the images $\bigoplus\limits_{\al=1}^s \im (\chi_\al)$.
In addition,
$V_{\is}\mapsto
V_{\{i_s\}}\times V_{\{i_s+i_{s-1}\}}\times\ld\times V_{\{n\}}.$
Thus we obtain a map
$$
\Fl_{\is}\hk \Fl_{\{i_s\}}\times \Fl_{\{i_s+i_{s-1}\}}\times\ld
\times \Fl_{\{n\}}.
$$
Because of the lemmas $(\ref{s=1})$ and $(\ref{embed})$ we obtain
the $\G$-equivariant isomorphism $\Fl_{\is}\simeq \sh_{\is}$.
\end{proof}

\section{Line bundles on $\sh_{\is}$}
Fix the set of the curves
$C_j\hk\sh_{\is},
j=0,\ld,n-1$,
$$C_j=\ov{\{\exp (ze_j)\cdot [v_{\is}],\ z\in\C\}}\simeq \Pro^1.$$
Let $\E$ be the line bundle on $\sh_{\is}$. We write
$\E=\Ob(a_1,\ld,a_n)=\Ob(A)$  ($a_i\in\Z$) if
$$\E|_{C_j}=\Ob(a_1+\ld+ a_{n-j}).$$
In \cite{mi3} was shown that the bundle on $\shn$ is uniquely determined
by the numbers $a_i$. The same statement is true for the
general $\sh_{\is}$. In the case of $\shn$ for any $\{a_i\}$ there exists a
bundle $\Ob(a_1,\ld,a_n)$. The general case is considered in the proposition
$(\ref{linebun})$.

\begin{lem}
\label{deg}
Consider the homomorphism $h_{\is}:\shn\to \sh_{\is}$. Recall
the subvariety $N_{n-1}(\{1,\ld,1\})\simeq \sh^{(n-2)}\times\Pro^1$ of $\shn$.
Denote by $L$ the following projective line:
\begin{equation}
L=[v_{\{\underbrace{{\scriptstyle 1,\ld,1}}_{n-2}\}}]\times\Pro^1\hk
N_{n-1}(\underbrace{1,\ld,1}_n)\hk \shn.
\end{equation}
If $i_s>1$ then $h_{\is}$ maps  $L$ to the point.
\end{lem}
\begin{proof}
From the results of \cite{mi3} follows that
\begin{equation}
\label{L}
L=\Slt^{(n)} \cdot [v(n-2)\T v(n-3)\T \ld \T v(0)\T v(1)],
\end{equation}
where $\Slt^{(n)}$ stands for the $\Slt$ acting on the $n$-th (last) factor
of the tensor power $F^{\T n}$.

Recall (see lemma $(\ref{glue})$, corollary $(\ref{his})$) that the
map $h_{\is}$ can be regarded as a part of the following commutative diagram
(the horizontal arrows come from the lemma $(\ref{embed})$):
$$
\begin{CD}
\shn @>>> \sh_{\{n\}}\times\sh_{\{n-1\}}\times\ld\times \sh_{\{ 1\}}\\
@V{h_{\is}}VV  @VV{P}V\\
\sh_{\is} @>>> \sh_{\{n\}}\times\sh_{\{n-i_1\}}\times\ld\times \sh_{\{i_s\}}
\end{CD}
$$
where $P$ is a map of "forgetting" of some factors. In the case
$i_s>1$ the last, $n$-th factor is one to "forget". Thus because of
the formula $(\ref{L})$ $h_{\is}$ maps $L$ to the point.
\end{proof}

\begin{prop}
\label{linebun}
Let $A=(a_1,\ld,a_n)\in\Z^n$ be of the type $\{j_1,\ld,j_{s_1}\}$.
Then the bundle $\Ob(A)$ on $\sh_{\is}$ really exists if and only if
$\{j_1,\ld,j_{s_1}\}\le {\is}$.
\end{prop}
\begin{proof}
First note that for any $A=(a_1,\ld,a_n), a_i>0$ of the type less or equal
to $\is$
there exists homomorphism $\imath_A:\sh_{\is}\to \Pro (M^A)$. In fact, if
$a_{i_1+\ld+i_\al}\ne a_{i_1+\ld+i_\al+1}$ for all $\al$, then $\imath_A$
is an embedding coming from the realization $\sh_{\is}=\sh_A$. If there
exists $\al$ with $a_{i_1+\ld+i_\al}=a_{i_1+\ld+i_\al+1}$, then $\imath_A$
is a composition of the above embedding and a homomorphisms from the lemma
$(\ref{glue})$. It was shown in \cite{mi3} that
$$\imath_A^* \Ob(1)=\Ob(a_1-1,\ld,a_n-1).$$
Thus for any $A$ with $\{j_1,\ld,j_{s_1}\}\le \is$ and $a_i\ge 0$ the
line bundle $\Ob(A)$ really exists on $\sh_{\is}$.
The case of an arbitrary $A$ is an immediate consequence, because any $A$
of the type $\{j_1,\ld,j_{s_1}\}$ can be represented as $A=B-C$
($a_i=b_i-c_i$), where $B,C\in\N^n$.

Now ewe need to prove that the bundle $\Ob(A)$ really exists on $\sh_{\is}$
only if the type of $A$ is less or equal to $\is$.
Recall that
$$\ov{\left\{\exp\left(\sum\nolimits_{j=1}^{n-1} z_je_j\right)\cdot
[v_{\is}],\ z_j\in\C\right\}}\simeq
\sh_{\{i_1,\ld,i_{s-1},i_s-1\}}.$$
The restriction of $\E=\Ob(a_1,\ld,a_n)$ to this subvariety equals to
$\Ob (a_1,\ld,a_{n-1})$.
Using the induction assumption we know that
$\{j_1,\ld,j_{s_1}-1\}\le\{i_1,\ld, i_s-1\}$.
Note that in the case $i_s=1$ we obtain that the type of $A$ is less or equal
to $\is$.

Let $i_s>1$. We want to prove that $j_{s_1}>1$.
Consider the variety $\shn$ and the line $L\hk\shn$
(see lemma $(\ref{deg})$).
It was shown in \cite{mi3} that the restriction of the bundle
$\Ob (A)$ on $\shn$ to this line equals to
$\Ob(a_n-a_{n-1})$. Using the lemma $(\ref{deg})$ we obtain that
if $\Ob(A)$ really exists on $\sh_{\is}$ then $a_n-a_{n-1}=0$.
Thus $j_{s_1}>1$.
Theorem is completely proved.
\end{proof}

\begin{cor}
\label{fus}
Let $A\in\N^n$ be of the type less or equal to $\is$.
Then
$$H^0(\Ob(a_1-1,\ld,a_n-1),\sh_{\is})\simeq (M^A)^*.$$
\end{cor}
\begin{proof}
It was proved in \cite{mi3} that
$$H^0(\Ob(a_1-1,\ld,a_n-1),\shn)\simeq (M^A)^*.$$
The statement of the corollary follows from the fact that
$$h_{\is}^* \Ob(a_1-1,\ld,a_n-1)=\Ob(a_1-1,\ld,a_n-1)$$ and
the induced map of the sections is an isomorphism.
\end{proof}

\section{Infinite-dimensional consequences}
Recall that in \cite{mi3} for any $A=(a_1\le \ld\le a_n)\in\N^n$ the
infinite-dimensional variety $Gr_A$ was defined as an inductive limit of
the Schubert varieties. Namely, let
$$A^{(i)}=(a_1,\ld,a_n,\underbrace{a_n,\ld,a_n}_{2i}).$$
Then we have an embedding $\sh_{(A^{(i)})}\hk \sh_{(A^{(i+1)})}$ and
$Gr_A=\lim_{i\to\infty} \sh_{(A^{(i)})}$.
Because of the theorem $(\ref{isom})$ we have a collection of varieties
$Gr_{\is}$:
$$Gr_{\is}=\lim_{i\to\infty} \sh_{\{i_1,\ld,i_s+2i\}}.$$
By definition we have an action of the group $\Slt(\C[t])$ on
$Gr_{\is}$. It was shown in \cite{mi2, mi3} that in fact we have an
action of the group $\Slth$ (the reason is that the space
$\lim_{i\to\infty} M^{A^{(i)}}$ is not only $\sltc$ module, but also
has a structure of a representation of the Lie algebra $\slth$).

\begin{lem}
Let $\jmath: \sh_{\is}\hk \sh_{\{i_1,\ld,i_s+2\}}$,
$B=(b_1,\ld,b_n)\in\Z^n$ the element of the type $\is$.
Then
$\jmath^* \Ob(B^{(1)})=\Ob(B).$
\end{lem}
\begin{proof}
It is enough to prove our lemma for $B\in (\N\setminus 1)^n$.
In this case
the embedding $\jmath$ is the restriction of the embedding of the
projective spaces
$\widehat\jmath: \Pro (M^B)\hk \Pro (M^{B^{(1)}})$, constructed
in \cite{mi3} (we use the embeddings $\imath_B: \sh_{\is}\hk \Pro(M^B)$ and
$\imath_{B^{(1)}}:\sh_{\{i_1,\ld,i_s+2\}}\hk \Pro (M^{B^{(1)}})$).
Thus
$\imath_{B^{(1)}}\jmath=\widehat\jmath \imath_B$.
We obtain
\begin{equation*}
\jmath^* \Ob(B^{(1)})=
\jmath^* \imath_{B^{(1)}}^* \Ob_{\Pro (M^{B^{(1)}})}(1)=
\imath_B^* \widehat\jmath^* \Ob_{\Pro (M^{B^{(1)}})}(1)=
\imath_B^* \Ob_{\Pro (M^B)}(1)=\Ob(B).
\end{equation*}
Lemma is proved.
\end{proof}

Now let ${\bf E}$ be a line
bundle on $Gr_{\is}$. We write
$${\bf E}=\Ob(B^{(\infty)})\ \text{ if } \
{\bf E}|_{\sh_{\{i_1,\ld,i_s+2i\}}}=
\Ob(B^{(i)}).$$
Let $B\in\N^n$.
Consider the projective limit of the spaces of sections
\begin{multline*}
H^0(\Ob(B^{(\infty)}),Gr_{\is})=\varprojlim_{i}
H^0(\Ob(B^{(i)}),\sh_{\{i_1,\ld,i_s+2i\}})
=\\ =\varprojlim_{i}
\left(M^{(b_1+1,\ld,b_n+1,(b_n+1)^{2i})}\right)^*.
\end{multline*}
It was proved in \cite{mi2} that we have an isomorphism of $\slth$ modules:
\begin{equation*}
\varinjlim_i M^{(b_1+1,\ld,b_n+1,(b_n+1)^{2i})}\simeq
\bigoplus_{j=0}^{b_n} c_{j;b_1,\ld,b_n} L_{j,b_n},
\end{equation*}
where $L_{j,b_n}$ are level $b_n$ irreducible $\slth$ modules and
the numbers $c_{j;b_1,\ld,b_n}$ are defined by the following equation
in the Verlinde algebra of the level $b_n+1$, associated with the Lie algebra
$\slt$:
$$[b_1]\cdot\ld \cdot [b_n]=\sum_{j=0}^{b_n} c_{j;b_1,\ld,b_n}[j]$$
(here the notation $[j]$ stands for the element of the Verlinde algebra,
corresponding to the $[j+1]$-dimensional irreducible representation of
$\slt$).
Thus we obtain the following proposition:
\begin{prop}
\label{inf}
Let $B\in\N^n$. Then we have an isomorphism of $\slth$ modules
$$H^0(\Ob(B^{(\infty)}),Gr_{\is})\simeq
\bigoplus_{j=0}^{b_n} c_{j;b_1,\ld,b_n} L_{j,b_n}^*.$$
\end{prop}

\section {Discussion of the singularities of the Schubert varieties}
Recall that in \cite{mi3} was shown that $\shn$ is a smooth variety.
\begin{lem}
Let $i_1+\ld +i_s=n$ and $s\ne n$. Then $\sh_{\is}$ is a singular variety.
\end{lem}
\begin{proof}
Recall (see \cite{mi3}) that there is a bundle
$$\widetilde\pi:\G\cdot [v_{\is}]\to \Slt\cdot [v_{\is}]\simeq \Pro^1$$
with a fiber $\C^{n-1}$. It is easy to show that the closure of the fiber
$\widetilde\pi(x)^{-1}$ is isomorphic to $\sh_{\{i_1,\ld,i_{s-1}, i_s-1\}}$
for any $x\in\Pro^1$. Note that because $\widetilde\pi$ is
$\Slt$-equivariant homomorphism, $\sh_{\is}$ is smooth if and only if
\begin{enumerate}
\item $\ov{\widetilde\pi(x)^{-1}}$ is smooth for any $x\in\Pro^1$,
\item $\ov{\widetilde\pi(x)^{-1}}\cap \ov{\widetilde\pi(y)^{-1}}=\emptyset$
if $x\ne y$.
\end{enumerate}
Suppose that we have already proved our lemma for $\sum i_\al <n$. Then the
first condition means that $i_1=\ld=i_{s-1}=1, i_s\le 2$. In \cite{mi3} was
proved that the condition $(2)$ holds for $i_s=1$. In the same way one can
prove that otherwise the closures of the fibers of $\widetilde\pi$ intersect.
Thus $i_s=1$ and lemma is proved.
\end{proof}

It is interesting to describe the variety of the singular points of
$\sh_{\is}$ in the general case. In the next proposition
we consider the case of $\sh_{\{n\}}$.
\begin{prop}
Recall (see proposition $(\ref{mainprop}))$ that
$\sh_{\{n\}}=\G\cdot [v_{\{n\}}]\sqcup N$ and $N\simeq \sh_{\{n-2\}}$.
\begin{enumerate}
\item $N$ is a variety of the singular points of $\sh_{\{n\}}$.
\item There exists a line bundle $\E$ on $\sh_{\{n\}}$ and
$s_1,s_2\in H^0 (\E, \sh_{\{n\}})$ such that $N=\{s_1=0\}\cap \{s_2=0\}$.
\end{enumerate}
\end{prop}
\begin{proof}
For the proof of the first part of the proposition it is enough to show
that for
any $x\in\Pro^1$  variety $N$ is contained in the closure
$\ov{\widetilde\pi^{-1}(x)}$.  First let us show the latter for the point
$x=[v_{\{n\}}]$ (recall that we have a realization
$\Pro^1=\Slt\cdot [v_{\{n\}}]$).
Note that the case of an arbitrary $x$ is an immediate consequence, because
$\widetilde\pi$ is an $\Slt$-equivariant map.

Recall that $\sh_{\{n\}}= \sh_{(2^n)}$.
From the proof of the
proposition $(\ref{mainprop})$ in the case of $\sh_{\{n\}}$ follows that
\begin{equation*}
\lim_{z\to\infty} \exp(e_{n-1}z) [v_{\{n\}}]=[e_{n-1} v_{(2^n)}]=
[v_{\{n-2\}}]\in \sh_{\{n-2\}}\simeq N
\end{equation*}
(recall that $e^2_{n-1} v_{(2^n)}=0$). Thus because of
$$N=\ov {\Slt(\C[t]/t^{n-2})\cdot [e_{n-1} v_{(2^n)}]}$$
we obtain $N\hk \ov{\widetilde\pi^{-1}([v_{\{n\}}])}$. The first part of the
proposition is proved.

Recall that $\Pic(\sh_{\{n\}})\simeq \Z$. Let $\E=\Ob(1,\ld,1)$ be the
generator of this group. We have an isomorphism (see corollary
$(\ref{fus})$):
\begin{equation}
\label{section}
H^0(\E,\sh_{\{n\}})\simeq \left(M^{(2^n)}\right)^*.
\end{equation}
Fix some basis $\{b_i\}_{i=1}^{2^n}$ of $M^{(2^n)}$ which is homogeneous with
respect to the $h_0$-grading and $b_1=v_{(2^n)}, b_{2^n}=u_{(2^n)}$.
Let $s_1,s_2\in H^0(\E,\sh_{\{n\}})$, $s_1=v_{\{n\}}^*, s_2=u_{\{n\}}^*$
(thus $s_1, s_2$ are the elements of the dual basis $b_i^*$).
Then
\begin{equation}
\label{annih}
\{s_1=0\}=\ov{\widetilde\pi^{-1}([u_{\{n\}}])},\quad
\{s_2=0\}=\ov{\widetilde\pi^{-1}([v_{\{n\}}])}.
\end{equation}
Let us prove that the first equality really holds (the proof for the second
one is quite similar).
In fact,
\begin{equation}
\label{union}
\sh_{\{n\}}=\left\{\exp\left(\sum_{j=0}^{n-1} e_jz_j\right)\cdot [v_{\{n\}}],
z_j\in\C\right\}\bigsqcup
\ov{\left\{\exp\left(\sum_{j=1}^{n-1} z_jf_j\right)\cdot
[u_{\{n\}}], z_j\in\C\right\}}.
\end{equation}
Rewrite $(\ref{union})$ as $\sh_{\{n\}}=R_1\sqcup R_2$.
Recall the embedding $\imath_{(2^n)}:\sh_{\{n\}}\hk \Pro (M^{(2^n)})$ and
the equality $\imath_{(2^n)}^* \Ob(1)=\E$. Note that
$(\ref{section})$ means that the restriction map
$H^0(\Ob(1), \Pro (M^{(2^n)}))\to H^0 (\E, \sh_{\{n\}})$ is an
isomorphism. Thus
$\{s_1=0\}\cap R_1=\emptyset$. In addition
$R_2\hk \Pro (\bra b_i, i>1\ket)$ and so $s_1|_{R_2}=0$.
But $R_2=\ov{\widetilde\pi^{-1}([u_{\{n\}}])}$.
We have proved
$(\ref{annih})$ and hence our proposition is completely proved.
\end{proof}

\end{document}